\documentclass[12pt,a4paper]{article}
\usepackage{amsmath,amssymb,amsfonts}
\def\Bbb{\mathbb}

\title{\bf The Petersson-Knopp identity and Farey neighbours}

\author{Kurt Girstmair}
\date{}

\makeatletter
\let\@@maketitle=\maketitle
\def\maketitle{\def\thispagestyle##1{\relax}\@@maketitle}
\makeatother
\newtheorem{theorem}{Theorem}

\newtheorem{lemma}{Lemma}

\def\BE{\begin{equation}}
\def\EE{\end{equation}}
\def\BD{\begin{displaymath}}
\def\ED{\end{displaymath}}
\def\BA{\begin{array}}
\def\EA{\end{array}}
\def\BEA{\begin{eqnarray}}
\def\EEA{\end{eqnarray}}
\def\BI{\bibitem}

\def\Z{\Bbb Z}

\def\R{\Bbb R}

\def\phi{\varphi}

\def\MB{\mbox}
\def\LD{\ldots}
\def\OV{\overline}

\def\DIV{\,|\,}
\def\NDIV{\, \nmid \,}

\def\BQ{``}
\def\EQ{'' }

\def\NI{\noindent}
\def\MN{\medskip\noindent}
\def\STOP{\hfill$\Box$}

\def\DED{Dedekind }

\begin{document}
\maketitle

\begin{abstract}
\noindent
We study Dedekind sums $S(a,b)$ near Farey points of the interval $[0,b]$. Each of these Dedekind sums is connected with a
set of other Dedekind sums by the Petersson-Knopp identity. In the case considered here, this identity has a very specific interpretation,
inasmuch as each Dedekind occurring in this identity is close to a certain expected value. Conversely, each of these expected values occurs with a certain frequency,
a frequency that is consistent with the Petersson-Knopp identity.
\end{abstract}

\section*{1. Introduction}

Let $b$ be a positive integer and $a\in \Z$. The classical \DED sum $s(a,b)$ is defined by
\BD
   s(a,b)=\sum_{k=1}^b ((k/b))((ak/b))
\ED
where $((\LD))$ is the \BQ sawtooth function\EQ defined by
\BD
\label{1.1}
  ((t))=\begin{cases}
                 t-\lfloor t\rfloor-1/2, & \MB{ if } t\in\R\smallsetminus \Z; \\
                 0, & \MB{ if } t\in \Z.
               \end{cases}
\ED
(see, for instance, \cite{RaWh}). In many cases it is more
convenient to work with
\BD
 S(a,b)=12s(a,b)
\ED instead.
We call $S(a,b)$ a {\em normalized} \DED sum. In addition, we say that $S(a,b)$ a {\em primitive} \DED sum,
if $(a,b)=1$. In the opposite case $S(a,b)$ is called {\em imprimitive}. Since
\BD
 S(ad,bd)=S(a,b)
\ED
for every positive integer $d$ (see \cite[Th. 1]{RaWh}), each imprimitive \DED sum  $S(a,b)$ is equal to the primitive
\DED sum $S(a/d,b/d)$, where $d=(a,b)$. We also note the periodicity
\BE
\label{1.1.0}
   S(a+b,b)=S(a,b)
\EE
of (not necessarily primitive) \DED sums.

Let us start with a special case of what we are doing in the sequel. Let $a<b$ be positive integers, $(a,b)=1$, and $p$ a prime not dividing $a, b$.
Then the normalized \DED sums
\BE
\label{1.1.1}
 S(pa,b)\MB{ and } S(a+jb,pb), j\in\{0,\LD,p-1\},
\EE
are primitive up to one exception. Indeed, if $a+jb\equiv 0\mod p$, then $S(a+jb,pb)=S((a+jb)/p,b)$. Suppose we know that all \DED sums (\ref{1.1.1}) are
positive. Then we also know that $S(a,b)$ is positive. Moreover, we know that at least one of the \DED sums (\ref{1.1.1})
is $\ge S(a,b)$, whereas the sum of any $p$ of them must be $<(p+1)S(a,b)$. This is an immediate consequence of the Petersson-Knopp identity, which, in this
special case, reads
\BD
  S(pa,b)+\sum_{j=0}^{p-1}S(a+jb,pb)=(p+1)S(a,b).
\ED
In what follows we discuss a situation where we know much more, namely, that one of the \DED sums (\ref{1.1.1}) is close to $pS(a,b)$, whereas each of the $p$ remaining ones is
close to $S(a,b)/p$. Hence the Petersson-Knopp identity has a very specific interpretation in this context.

In two previous papers \cite{Gi1,Gi2} we studied the behaviour of primitive \DED sums near {Farey points}. We briefly recall the
necessary notation. Let the positive integer $b$ be given and assume $b\ge 4$. For a positive integer $d$, $d<b^{1/3}$, let $c\in\Z$, $(c,d)=1$.
Then $c/d$ is a Farey fraction of an order $<b^{1/3}$ in the usual sense (see, \cite[p. 125]{Hu}).
We say that $b\cdot c/d$ is a {\em Farey point with respect to} $b$.
Put
\BE
\label{1.2}
  \alpha=\sqrt{b/d^3}.
\EE
We call
\BD
\label{1.3}
   \{x\in\R: |x-b\cdot c/d|\le \alpha-1\}
\ED
the Farey interval belonging to $b\cdot c/d$.
Now let $a$ be an integer, $(a,b)=1$, inside the Farey interval. Then the primitive \DED sum $S(a,b)$ is $<0$, if $a<b\cdot c/d$, and $S(a,b)>0$, if $a>b\cdot c/d$ (see \cite[Th. 1 and formula (5)]{Gi1}).
In order to avoid tedious distinctions, we restrict ourselves to integers $a$ in the {\em right} half of the Farey interval, so $S(a,b)>0$. The whole theory remains
valid for integers in the left half, but with $S(a,b)$ negative.

Hence we say that $a\in\Z$, $(a,b)=1$, is a {\em Farey neighbour} of the point $b\cdot c/d$ if
\BE
\label{1.3.1}
0\le a-b\cdot c/d\le \alpha-1.
\EE
Note that $a-b\cdot c/d\ne 0$ since $a/b=c/d$ is impossible (both fractions are reduced, and $0<d<b$).
For a Farey neighbour $a$, $S(a,b)$ is not only positive, but its value is,
as a rule, close to an expected value, which can be defined as follows. Put
\BE
\label{1.4}
  q=ad-bc.
\EE
Then $q>0$ since $q/d=a-b\cdot c/d>0$. Now the {\em expected value} of $S(a,b)$
is
\BE
\label{1.5}
 E(a,b)= \frac b{dq}
\EE
(which is $>0$).
In Section 3 we will see why $S(a,b)$ is, in general, close to $E(a,b)$ if $a$ is a Farey neighbour of $b\cdot c/d$.

The Petersson-Knopp identity (see \cite{Pa}) is a relation between $S(a,b)$ and certain other \DED sums. Indeed, if $n$ is a natural number,
then
\BE
\label{1.7}
   \sum_{r\DIV n}\sum_{j=0}^{r-1} S\left(\frac nr a+jb,rb\right)=\sigma(n)S(a,b).
\EE
Here $r$ runs through the (positive) divisors of $n$ and $\sigma(n)=\sum_{r\DIV n}r$ is the sum of the divisors of $n$.

The \DED sums in (\ref{1.7}) are not necessarily primitive. In order to apply results about Farey neighbours, we need primitive \DED sums, however.
In view of the periodicity (\ref{1.1.0}), it suffices to restrict $c$ to the range
$0\le c<d$, $(c,d)=1$. Let $a$ be a Farey neighbour of $b\cdot c/d$. For $r\DIV n$ and $j\in\{0,\LD,r-1\}$ put
\BD
\label{1.9}
  k(r,j)=\left(\frac nr a+jb, rb\right) \MB{ and } m(r,j)=\left(\frac nr c+jd,rd\right).
\ED
So both $k(r,j)$ and $m(r,j)$ are positive integers.
Moreover, put
\BD
\label{1.11}
 a(r,j)=\frac{\frac nr a+jb}{k(r,j)},\, b(r,j)=\frac{rb}{k(r,j)},\, c(r,j)=\frac{\frac nr c+jd}{m(r,j)},\, d(r,j)=\frac{rd}{m(r,j)}.
\ED
In the sequel we simply write
\BD
S[r,j]=S(a(r,j),b(r,j))=S\left(\frac nra+jb,rb\right)
\ED
and
\BD
 E[r,j]=E(a(r,j), b(r,j)).
\ED
Then we have the following result:

\begin{theorem} 
\label{t1}
In the above setting, let $0\le c<d$, $(c,d)=1$, $\alpha\ge n^{3/2}+n$ and
\BD
  0<a-b\cdot c/d\le\alpha/n-1.
\ED
For each pair $(r,j)$, $r\DIV n$, $j\in\{0,\LD,r-1\}$, the number $a(r,j)$ is a Farey neighbour of $b(r,j)\cdot c(r,j)/d(r,j)$. Hence
$S[r,j]$ is positive. Its expected value is
\BE
\label{1.13}
        E[r,j]=\frac{m(r,j)^2}n \cdot E(a,b),
\EE
where $E(a,b)$ is the expected value of $S(a,b)$, see {\rm(\ref{1.5})}.
\end{theorem} 

In view of the Petersson-Knopp identity (\ref{1.7}), one expects that
\BE
\label{1.14}
  \sum_{r\DIV n}\sum_{j=0}^{r-1}E[r,j]=\sigma(n)E(a,b).
\EE
This is true, but we have a much more precise result about the expected values $E[r,j]$. Indeed, they follow
a very regular pattern.

\begin{theorem} 
\label{t2}
In the above setting, the numbers $m(r,j)$ divide $n$. Conversely, for every positive divisor $m$ of $n$,
\BE
\label{1.15}
  \#\left\{(r,j); r\DIV n, j\in\{0,\LD,r-1\}, E[r,j]=\frac{m^2}n E(a,b)\right\}=\frac nm.
\EE
\end{theorem} 

By (\ref{1.13}) and (\ref{1.15}), the left hand side of (\ref{1.14}) reads
\BD
 \sum_{m\DIV n} \frac nm\cdot\frac{m^2}nE(a,b),
\ED
which obviously equals $\sigma(n)\cdot E(a,b)$.

\MN
{\em Example.} Let $n=12$. In this case there are $\sigma(12)=28$ \DED sums $S[r,j]$.
The corresponding values of $E[r,j]/E(a,b)$ are $1/12, 1/3, 3/4, 4/3, 3, 12$, respectively.
Let $d=9$ and $c=1$.
We choose $b$ so large that $\alpha/n-1\ge 10$. This means $b\ge 12702096$. Then it is obvious that
$\alpha\ge 132\ge n^{3/2}+n\approx 53.569$. We have used a random generator to produce a number $b$, $1.28\cdot 10^7<b<10^8$. It has given us
$b=31537789$. The Farey point $b\cdot c/d$ is approximately $3504198.78$. Since $\alpha/n-1\approx 16.33$, we can choose $a=3504214$,
which is prime to $b$, and $a-b\cdot c/d\approx 15.22$.  Then $S(a,b)\approx 25537.432$ and $E(a,b)\approx 25578.093$. We have computed the
relative deviation
\BE
\label{1.17}
  \left|\frac{S[r,j]}{E[r,j]}-1\right|
\EE
of each of the said 28 \DED sums from its expected value. It turns out that the largest relative deviation is $\approx 0.04659$ or nearly $4.7$ percent.
It occurs for $r=6$, $j=1$, where $E[r,j]=(1/12)E(a,b)$. The mean relative deviation, i.e., the arithmetic mean of all values (\ref{1.17}),
is $\approx 0.0060$ or $0.6$ percent. Further empirical results can be found in Section 3.

\MN
{\em Remark}. The example shows that there are, compared with the size of $b$, only few integers $a$ such that $0<a-b\cdot c/d\le \alpha/n-1$
for a fixed value of $d$ and $0\le c<d$, $(c,d)=1$. In the case of the example their number amounts to $\approx 6\cdot 15=90$.
However, one should be aware of the fact that
each number $a$ of this kind also satisfies $a-b\cdot c/d\le\alpha/n'-1$ for all integers $n'$, $1\le n'<n$. Therefore, if $(a,b)=1$,
the number $a$ gives rise not only to the $\sigma(n)$ \DED sums $S[r,j]$ for $n$,
but also to $\sigma(n')$ analogous \DED sums for each positive integer $n'<n$ (the case $n'=1$ includes $S(a,b)$).
For $n=12$ their totality amounts to $\sigma(1)+\sigma(2)+\LD+\sigma(12)=112$.
In general,
\BD
  \sum_{n'=1}^n\sigma(n') = \frac{\pi^2}{12}n^2+O(n\log n),
\ED
see \cite[p. 113]{Hu}.
Hence there is quite a number of \DED sums whose expected values are known.

\section*{2. Proofs}

Let the assumptions of Theorem \ref{t1} hold. In particular, let $r$ divide $n$ and $j\in\{0,\LD,r-1\}$.

We first show that $m(r,j)$ divides $n$. Let $p$ be a prime. We use the {\em $p$-exponent} $v_p(t)$ of an integer $t\ne 0$,
which is given by $t=p^{v_p(t)}t'$, $(p,t')=1$. We show that $v_p(m(r,j))\le v_p(n)$ for all primes $p$. To this end
recall that $m(r,j)=(\frac nr c+jd, rd)$.
First suppose $p\NDIV d$. Then $v_p(m(r,j))\le v_p(r)\le v_p(n)$. Next let $p\DIV d$, so $v_p(d)=s\ge 1$.
Since $(c,d)=1$, $v_p(c)=0$ and $v_p(\frac nr c)=v_p(\frac nr)$. If $v_p(\frac nr)<s$, then $v_p(\frac nrc+jd)=v_p(\frac nr)\le v_p(n)$. If $v_p(\frac nr)\ge s$,
then $v_p(n)\ge v_p(r)+s$. In this case $v_p(rd)=v_p(r)+s\le v_p(n)$, and $v_p(m(r,j))\le v_p(rd)\le v_p(n)$.

The same arguments work for $k(r,j)=(\frac nr a+jb, rb)$ and $a,b$ instead of $(c,d)$. They show that $k(r,j)$ divides $n$.

\MN
{\em Proof of Theorem \ref{t1}.} In order to simplify the notation for the purpose of this proof,
we write $a'=a(r,j)$, $b'=b(r,j)$, $c'=c(r,j)$, $d'=d(r,j), k'=k(r,j)$ and $m'=m(r,j)$.
First we observe $b'\ge b/k'\ge b/n$, and since $\alpha\ge 2n$, we have $b'\ge 4$.

Next we consider
\BD
\label{2.1}
  q'=a'd'-b'c'.
\ED
A short calculation shows
\BE
\label{2.3}
 q'=\frac n{k'm'}\,q,
\EE
where $q=ad-bc$, see (\ref{1.4}).
Now $a'$ is a Farey neighbour of $b'\cdot c'/d'$,
if $0<a'-b'\cdot c'/d' \le\sqrt{b'/d'^3}-1$, i.e.,
\BD
  0<q'\le \sqrt{{b'}/{d'}}-d',
\ED
see (\ref{1.3.1}). Here $q'>0$ follows from (\ref{2.3}), since $q>0$. Because $\sqrt{b'/d'}=\sqrt{m/k}\cdot\sqrt{b/d}$, $a'$ is a Farey neighbour of $b'\cdot c'/d'$, if
\BD
   \frac n{k'm'}\,q\le \sqrt{\frac{m'}{k'}}\cdot\sqrt{\frac bd}-\frac{rd}{m'}
\ED
by (\ref{2.3}). This condition can be written
\BE
\label{2.5}
  a-b\cdot\frac cd= \frac qd\le \frac{k'^{1/2}m'^{3/2}}n\cdot\alpha-\frac{rk'}n.
\EE
Let $\rho$ be the right hand side of (\ref{2.5}), i.e.,
\BD
 \rho=\frac{k'^{1/2}m'^{3/2}}n\cdot\alpha-\frac{rk'}n.
\ED
If $k'=m'=1$ and $r=n$, then $\rho$ becomes $\alpha/n-1$.
We show that $\rho$ is always $\ge \alpha/n-1$, provided that $\alpha\ge n^{3/2}+n$. In this case the condition
$q/d\le \alpha/n-1$ implies that $a'$ is a Farey neighbour of $b'\cdot c'/d'$ for {\em all} $r,j$ in question.

In the case $k'=1$ we have $rk'/n\le r/n\le 1$ and $\rho\ge \alpha/n-1$.
Hence assume $k'>1$. Since $\rho$ is $\ge k'^{1/2}\alpha/n-rk'/n$, $\rho<\alpha/n-1$ implies
$k'^{1/2}\alpha-rk'/n<\alpha/n-1$. Because $k'^{1/2}>1$, this inequality can be written
$\alpha<(rk'-n)/(k'^{1/2}-1)$. Since $r\le n$, it implies $\alpha<n(k'^{1/2}+1)$. We know that $k'$ divides $n$,
hence we obtain $\alpha<n^{3/2}+n$ as a necessary condition for $\rho<\alpha/n-1$.

Finally, we compute
\BD
  E(a',b')=\frac{b'}{d'q'}=\frac{rb/k'}{rd/m'\cdot q\cdot n/(k'm')}=\frac{m'^2}n\, \frac b{dq}=\frac{m'^2}nE(a,b).
\ED
\STOP

In the sequel we need the following notation. For positive integers $r$ and $d$ let $(r)_d$ and $(r)_d^{\bot}$ denote the {\em $d$-part} and the {\em $d$-free part} of $r$, respectively, i.e.,
\BD
   (r)_d=\prod_{p\DIV r,\, p\DIV d} p^{v_p(r)}\MB{ and }(r)_d^{\bot}=\prod_{p\DIV r,\, p\NDIV d} p^{v_p(r)},
\ED
where $v_p(r)$ is defined as above.
The proof of Theorem \ref{t2} is more complicated than that of Theorem \ref{t1} and based on the following lemmas.

\begin{lemma} 
\label{l1}
Let $r,d$ be positive integers and $s\in\Z$ such that $(s,d)=1$.
Then
\BD
  \#\{\OV k\in\Z/r\Z: \OV{s+kd}\in (\Z/r\Z)^{\times}\}=(r)_d\,\phi((r)_d^{\bot}),
\ED
where $\phi$ denotes Euler's totient function.
\end{lemma} 

\NI
{\em Proof.} We use the Chinese remainder theorem to decompose $\Z/r\Z$ into its $p$-parts $\Z/p^{e_p}\Z$, where $e_p=v_p(r)\ge 1$.

Case 1: $p\DIV d$. Then we have, for all $k\in \Z$, $s+kd\equiv s\not\equiv 0\mod p$, i.e.,
$\OV{s+kd}\in(\Z/p^{e_p})^{\times}$. Hence
\BD
  \#\{\OV k\in\Z/p^{e_p}\Z: \OV{s+kd}\in(\Z/p^{e_p}\Z)^{\times}\}=p^{e_p}.
\ED

Case 2: $p\NDIV d$. Let $k\in\Z$. Let $d^*$ be an inverse of $d$ mod $p$. Then $s+kd\not\equiv 0\mod p$, if, and only if, $k\not\equiv -sd^*\mod p$. Therefore,
\BD
  \#\{\OV k\in\Z/p^{e_p}\Z: \OV{s+kd}\in(\Z/p^{e_p}\Z)^{\times}\}=p^{e_p}(1-1/p)=\phi(p^{e_p}).
\ED
\STOP

\begin{lemma} 
\label{l2}
Let $n$ be a positive integer and $m>0$ a divisor of $n$. Let $c,d\in\Z$, $0\le c<d$, $(c,d)=1$, and $\delta=(m,d)$.
Put $n'=n/\delta$, $m'=m/\delta$ and $d'=d/\delta$. Then
\BE
\label{2.7}
  \#\left\{(r,j): r\DIV n,0\le j\le r-1, m=\left(\frac nr c+jd,rd\right)\right\}=
  \sum_{\genfrac{}{}{0pt}{1}{m'\DIV r\DIV n'}{(n/r,d)=\delta}}(r/m')_{d'}\,\phi((r/m')_{d'}^{\bot}).
\EE
\end{lemma} 

\NI
{\em Proof.} We determine, for given positive divisors $m, r$ of $n$,
\BE
\label{2.9}
  \#\left\{j:0\le j\le r-1, m=\left(\frac nr c+jd,rd\right)\right\}.
\EE
First we show that (\ref{2.9}) equals $0$ if $(n/r,d)\ne\delta$. To this end suppose that $m=(\frac nr c+jd,rd)$ for some $j$.
Since $\delta\DIV m$, we have $\delta\DIV \frac nr c+jd$, and because $\delta\DIV d$, we obtain $\delta\DIV \frac nrc$.
But $(c,d)=1$, and so $\delta\DIV \frac nr$. Put $d_r=(\frac nr,d)$. We have seen $\delta\DIV d_r$. Conversely, $d_r$ divides
both $\frac nr c+jd$ and $rd$, whence $d_r\DIV m$. But $d_r\DIV d$, which implies $d_r\DIV (m,d)=\delta$. Altogether, $d_r=\delta$.
This means that $m=(\frac nr c+jd,rd)$ can hold only if $d_r=\delta$.

Therefore, we can restrict our investigation of (\ref{2.9}) to those $r$ for which $(\frac nr,d)=\delta$. As above,
put $d'=d/\delta$ and $n'=n/\delta$. Since $\delta\DIV n/r$, $r$ divides $n'$. Suppose that $m=(\frac nr c+jd,dr)$. Then
$ m=\delta m' \MB{ with } m'=(\frac{n'}r c+jd',rd')$.
Because $(\frac nr,d)=\delta$, we have $(\frac{n'}r,d')=1$ and $(\frac{n'}r c+jd',d')=1$. Accordingly,
\BE
\label{2.11}
   m'=\left(\frac{n'}r c+jd',r\right).
\EE
Conversely, suppose that $m'=m/\delta$ divides $r$. Since $(m',d')=1$, there is a number $j_0\in\{0,\LD,m'-1\}$ such that
$\frac{n'}r c+j_0d'\equiv 0\mod m'$. If $m'$ has the form (\ref{2.11}) for a number $j\in\{0,\LD,r-1\}$, then $j\equiv j_0\mod m'$, and so
$ j=j_0+km'$ for a uniquely determined $k\in\{0,\LD, r/{m'}-1\}$. For such a number $j$, we have
\BD
 \left(\frac{n'}r c+jd'\right)/{m'}=s+kd'
\ED
with $s=(\frac{n'}r c+j_0d')/m'$. Now (\ref{2.11}) holds if, and only if,
\BD
  (s+kd',r/{m'})=1.
\ED
Therefore, we have to count the $\OV k\in\Z/\frac r{m'}\Z$ such that $\OV{s+kd'}\in(\Z/\frac r{m'}\Z)^{\times}$. From Lemma \ref{l1}
we know that the number of these elements $\OV k$ equals
\BE
\label{2.13}
  \left(r/{m'}\right)_{d'}\,\phi\left(\left(r/{m'}\right)_{d'}^{\bot}\right).
\EE
This number equals that of (\ref{2.9}). We have to sum up the numbers (\ref{2.13}), observing that $(n/r,d)=\delta$. This yields
(\ref{2.7}).
\STOP

\MN
For positive integers $n, m$, $m\DIV n$, let $A(m,n)$ denote the number of (\ref{2.7}), i.e.,
\BD
  A(n,m)=\#\left\{(r,j): r\DIV n,0\le j\le r-1, m=\left(\frac nr c+jd,rd\right)\right\}.
\ED

\begin{lemma} 
\label{l3}
Let $n,m$ be positive integers, $m\DIV n$, and suppose $n=n_1n_2$ for positive integers $n_1$, $n_2$ such that $(n_1,n_2)=1$.
Put $m_1=(m,n_1)$ and $m_2=(m,n_2)$. Then
\BD
\label{2.15}
  A(n,m)=A(n_1,m_1)A(n_2,m_2).
\ED
\end{lemma} 

\NI
{\em Proof.} All entries of the right hand side of (\ref{2.7}) are multiplicative. Indeed, put $\delta_1=(\delta,n_1)$ and $\delta_2=(\delta,n_2)$. Then
$\delta=\delta_1\delta_2$. In the same way, $r=r_1r_2$ with $r_1=(r,n_1)$ and $r_2=(r,n_2)$. We also have $n'=n_1'n_2'$ with $n_1'=n_1/\delta_1$ and
$n_2'=n_2/\delta_2$. The respective identity holds for $m'$ and $m_1'=m_1/\delta_1$ and $m_2'=m_2/\delta_2$.
Further, $(n/r,d)=\delta$ if, and only if, $(n_1/r_1,d)=\delta_1$ and $(n_2/r_2,d)=\delta_2$. We note
$(r/m')_{d'}=(r_1/m_1')_{d_1'}(r_2/m_2')_{d_2'}$, where $d_1'=d/\delta_1$ and $d_2'=d/\delta_2$. The same identity holds when we apply the $\bot$ to the respective items.
Finally, the function $\phi$ is also multiplicative. In view of all that, we can write the sum over $r$ as the product of two sums over $r_1$ and $r_2$ and obtain the desired
result.
\STOP

\MN
{\em Proof of Theorem \ref{t2}.} We have to show that $A(m,n)=n/m$. By Lemma \ref{l3}, it suffices to prove this identity for prime powers $n=p^e$ and $m\DIV n$. Suppose that $m=p^k$,
$k\le e$, and $(d)_p=p^s$.

Case 1: $k\ge s$. Then $\delta=p^s$. We have $m'=p^{k-s}$ and $n'=p^{e-s}$. Let $r=p^t$ with $k-s\le t\le e-s$. By Lemma \ref{l2},
\BD
  A(n,m)=\sum_{\genfrac{}{}{0pt}{1}{k-s\le t\le e-s}{(p^{e-t},p^s)=p^s}}\phi(p^{t+s-k}),
\ED
since $r/m'=p^{t+s-k}$ and $(d')_p=(d/\delta)_p=1$. Obviously, $(p^{e-t},p^s)=p^s$ holds for all $t$ in question, because $e-t\ge s$. We obtain
\BD
   A(n,m)=\sum_{u=0}^{e-k}\phi(p^u)=p^{e-k}=n/m.
\ED

Case 2: $k<s$. Then $\delta=p^k$. Moreover, $m'=1$ and $n'=p^{e-k}$. If $r=p^t$, $0\le t\le e-k$, we have
\BD
  (n/r,d)=(p^{e-t}, p^s)=\begin{cases} p^{e-t}, \MB{ if } e-t\le s,\\
                                       p^s,     \MB{ if } e-t>s.
                          \end{cases}
\ED
Since $r$ must satisfy $(n/r,s)=\delta=p^k$ and $k<s$, only the first case is suitable for our purpose, and, indeed, only for $e-t=k$,
i.e., $t=e-k$. So only the summand for $r=p^{e-k}$ remains. We have $(d')_p=p^{s-k}$ with $s-k\ge 1$. Accordingly,
$(r/m')_{d'}=(r)_{d'}=r=p^{e-k}$ and $A(n,m)=n/m$, again.
\STOP

\section*{3. Theoretical and numerical evidence for the expected values}

It is a consequence of the three-term relation for $\DED$ sums that $S(a,b)$ is, in general, close to $E(a,b)=b/(dq)$ when $a$ is a Farey neighbour of $b\cdot c/d$.
Indeed, we have
\BE
\label{3.1}
 S(a,b)=\frac b{dq}+S(c,d)+S(t,q)+\frac d{bq}+\frac q{db}-3,
\EE
see \cite[Lemma 3]{Gi1}. Here $q$ is defined by (\ref{1.4}) and $t$ is an integer defined by $a,b,c,d$. The exact value of $t$ is not of interest for our purpose. First we observe
$d<b$ and, by (\ref{1.3.1}) and (\ref{1.4}), $q<\sqrt{b/d}<b$. We have, thus,
\BD
  0<\frac d{bq}+\frac q{db}<2.
\ED
Next we note
\BE
\label{3.3}
  |S(c,d)|<d \MB{ and } |S(t,q)|<q,
\EE
see \cite[Satz 2]{Ra}. In most cases, however, these \DED sums are much smaller, say $|S(c,d)|\le 5\log d$ and $|S(t,q)|\le 5\log q$.
Indeed, the main result of \cite{Va} allows to determine the asymptotic proportion of pairs $(c,d)$, $0\le c<d\le N$, $(c,d)=1$, such that
$|S(c,d)|<C\log d$ for a given constant $C>0$, as $N$ tends to infinity. For $C=5$ this proportion is about $76.8$ percent, and for $C=10$ about $88.0$ percent.

Another argument in favour of small values of $|S(c,d)|$ and $|S(t,q)|$ is the {\em mean value} of all \DED sums $S(c,d)$, $0\le c<d$, $(c,d)=1$, for a given
positive integer $d$. As $d$ tends to infinity, this mean value is $\le \log^2d\cdot 6/\pi^2+O(\log d)$, see \cite{GiSch}.

On the other hand, $b/(dq)> \sqrt{b/d}> b^{1/3}$, since $q<\sqrt{b/d}$ (recall $d<b^{1/3}$). These arguments support the hope that
the right hand side of (\ref{3.1}) is close to $E(a,b)=b/(dq)$ in most cases, a hope that is supported by empirical data, see below.

It should also be mentioned that the approximation of $E(a,b)$ becomes better when $d$ is small, say $d<b^{1/5}$, and the Farey neighbour $a$ tends to the Farey point $b\cdot c/d$.
Indeed, in this case $q/d$ tends to a positive value $\le 1$. So (\ref{3.3}) shows that the error caused by $S(c,d)$ and $S(t,q)$ has an absolute value $\le 2d<2b^{1/5}$.
On the other hand, $b/(dq)$ becomes $>s b/d^2>b^{3/5}$.

We return to the setting of the Theorems \ref{t1} and \ref{t2}.
Suppose that the size of $d$ is fixed, say $d\le n$, whereas $b$ may become large.
As in Theorem \ref{t1}, assume $\alpha\ge n^{3/2}+ n$ and $q/d\le \alpha/n-1$. Accordingly, all \DED sums
$S(a(r,j),b(r,j))=S[r,j]$ are positive for $r\DIV n$, $0\le j\le r-1$. The expected value of $S[r,j]$ equals $E[r,j]=(m(r,j)^2/n)E(a,b)$. By Theorem \ref{t2},
we know that $m(r,j)$ is a divisor of $n$, and, conversely,
each positive divisor $m$ of $n$ has the form $m=m(r,j)$ for exactly $n/m$ pairs $(r,j)$.

Empirical data shows that the relative deviation (\ref{1.17}) of $S[r,j]$ from $E[r,j]$ may be large, in the main, if $m(r,j)=k(r,j)=1$ and $q/d$ is close to $\alpha/n$. In this case $E[r,j]=(1/n)E(a,b)$.
This empirical observation can be explained as follows.
We have
\BD
  q(r,j)=a(r,j)d(r,j)-b(r,j)c(r,j)=\frac n{k(r,j)m(r,j)}\,q,
\ED
see (\ref{2.1}), (\ref{2.3}). Because $m(r,j)=k(r,j)=1$, we obtain $q(r,j)=nq$.
The influence of $S(c(r,j),d(r,j))$ on $S[r,j]$ in the sense of (\ref{3.1}) is limited since $|S(c(r,j),d(r,j))|\le d(r,j)\le rd\le n^2$. However, the influence of $S((t(r,j), q(r,j))$ may be large
if $q(r,j)=nq$ is close to $E[r,j]=(1/n)E(a,b)=b/(ndq)$, i.e., if $q/d$ is close to $\alpha/n$.

Let $(r,j)$ be of this kind and, in addition, the pair $(r_1,j_1)$ such that $m(r_1,j_1)\ge 2$. Then we have
\BD
  q(r_1,j_1)=\frac n{k(r_1,j_1)m(r_1,j_1)}\,q\le \frac{n}2\, q=\frac{q(r,j)}2.
\ED
On the other hand $E[r_1,j_1]\ge (4/n)E(a,b)=4E[r,j]$. This means
\BD
  \frac{q(r_1,j_1)}{E[r_1,j_1]}\le \frac 18 \cdot\frac{q(r,j)}{E[r,j]},
\ED
which is a much better proportion than $q(r,j)/E[r,j]$, in particular, in the bad case $q(r,j)\approx E[r,j]$.

As to empirical data, we have performed numerous computations, of which, however, we present only the case $n=12$ and $d=9$. We have computed the mean value of the relative
deviation (\ref{1.17}) both for all 28 pairs $(r,j)$, $r\DIV 12$, $j=0,\LD,r-1$, and only for those $(r,j)$ with $m(r,j)=1$ (and expected value $E[r,j]=(1/12)E(a,b)$). By the above, it is not surprising
that the first mean value is always smaller than the second.

We consider $b=10^8+k$, $1\le k\le 10000$, and choose the integer $a$ close to $b\cdot c/d+\alpha/n$. To be precise, $a$ is either $\lfloor b\cdot c/d+\alpha/n\rfloor-1$
or $\lfloor b\cdot c/d+\alpha/n\rfloor-2$. If none of these values of $a$ satisfies $(a,b)=1$, the number $b$ is ruled out. In this way there always remain $\ge 8000$ pairs $(b,a)$ to be investigated.
The following table lists the percentage of $b$'s such that the first mean value
\BD
   M_1=\frac 1{\sigma(n)}\sum_{r\DIV n}\sum_{j=0}^{r-1}\,\left|\frac{S[r,j]}{E[r,j]}-1\right|
\ED
is either $\ge 0.05$ or $<0.01$.
The table also displays the percentage of $b$'s such that the second mean value
\BD
   M_2=\frac 1{n}\sum_{m(r,j)=1}\,\left|\frac{S[r,j]}{E[r,j]}-1\right|
\ED
is either $\ge 0.1$ or $<0.01$.

\bigskip
\begin{tabular}{l||r|r|r|r|r|r} 
 c        &  1&  2&  4&  5&  7& 8\\
\hline
$M_1\ge 0.05$& 1.2 \%& 1.3 \%& 1.3 \%& 1.3 \%&1.3 \%&1.3 \% \rule{0mm}{5mm}\\
\hline
$M_1<0.01$& 93.4 \%& 93.7 \%& 93.7 \%& 93.6 \%& 93.8 \%& 92.9 \%\rule{0mm}{5mm}\\
\hline
$M_2\ge 0.1$& 1.2 \%& 1.3 \%& 1.3 \%& 1.3 \%&1.3 \%&1.3 \%\rule{0mm}{5mm}\\
\hline
$M_2<0.01$& 73.6 \%&80.5 \%&78.7 \%&78.8 \%& 80.6 \%& 70.6 \%\rule{0mm}{5mm}\\
\multicolumn{7}{c}{{\bf Table 1:} $n=12$, $d=5$, $b=10^8+k$, $1\le k\le 10000$\rule{0mm}{8mm}}
\end{tabular} 

\bigskip
We list the same data for numbers $b=10^9+k$, $1\le k\le 10000$.

\bigskip
\begin{tabular}{l||r|r|r|r|r|r} 
 c        &  1&  2&  4&  5&  7& 8\\
\hline
$M_1\ge 0.05$& 0.4 \%& 0.4 \%& 0.4 \%& 0.3 \%&0.3 \%&0.4 \%\rule{0mm}{5mm}\\
\hline
$M_1<0.01$& 97.9 \%&98.0 \%&98.0 \%&98.2 \%& 98.2 \%& 97.9 \%\rule{0mm}{5mm}\\
\hline
$M_2\ge 0.1$& 0.4 \%& 0.4 \%& 0.4 \%& 0.3 \%&0.3 \%&0.4 \%\rule{0mm}{5mm}\\
\hline
$M_2<0.01$& 94.4 \%&94.5 \%&94.7 \%&95.2 \%& 95.1 \%& 94.2 \%\rule{0mm}{5mm}\\
\multicolumn{7}{c}{{\bf Table 2:} $n=12$, $d=5$, $b=10^9+k$, $1\le k\le 10000$\rule{0mm}{8mm}}
\end{tabular} 

\bigskip
We obtain similar results when we use (pseudo-) random numbers $b$ of the same order of magnitude instead of the (more or less) consecutive numbers $b$ of the tables.
The tables suggest that the approximation of $E[r,j]$ by $S[r,j]$ becomes better when $b$ increases while $d$ and $n$ are fixed. This observation is supported by further computations.


\vspace{0.5cm}
\noindent
Kurt Girstmair            \\
Institut f\"ur Mathematik \\
Universit\"at Innsbruck   \\
Technikerstr. 13/7        \\
A-6020 Innsbruck, Austria \\
Kurt.Girstmair@uibk.ac.at

\end{document}